\documentclass{elsart}
\usepackage{amssymb}
\usepackage{amsmath}
\usepackage{amscd}
\usepackage[matrix,arrow]{xy}

\begin{document}

\begin{frontmatter}

\thanks[*]{The  author was supported in part by NSF Grants \#DMS-0200601 and 0500693}

\title{On the topology of the Kasparov groups and its applications}
\author[*]{Marius Dadarlat}
\address{Department of Mathematics,
 Purdue University, West Lafayette IN 47907, U.S.A.}
\ead{mdd@math.purdue.edu}
\ead[url]{http://www.math.purdue.edu/\~{}mdd}
\date{\today}
\begin{abstract} In this paper we establish a direct connection
 between stable approximate
unitary equivalence for $*$-homomorphisms and the topology of the KK-groups which avoids
entirely  C*-algebra extension theory and does not require
nuclearity assumptions.
 To this purpose we show that
  a topology on the
Kasparov groups  can be defined
 in terms of approximate
unitary equivalence for Cuntz pairs and that this topology coincides with both
Pimsner's topology and  the Brown-Salinas topology. We study the generalized
R{\o}rdam group $KL(A,B)=KK(A,B)/\bar{0}$, and  prove that if a separable exact
residually finite dimensional C*-algebra satisfies the universal coefficient
theorem in KK-theory, then it embeds in the UHF algebra of type $2^\infty$.
 In particular such an embedding exists for the C*-algebra of a second
countable amenable locally compact maximally almost periodic
group.
\end{abstract}

\begin{keyword} KK-theory\sep C*-algebras\sep amenable groups
\end{keyword}

\newcommand{\RG}{\ensuremath{R(\chi)}}
\newcommand{\MKB}{M(\mathcal{K}\otimes B)}
\newcommand{\bn}{\bar{\mathbb{N}}}
\newcommand{\SF}[1]{\underset{\mathcal{F}\!\!,#1}{\sim}}
\newcommand{\SFn}[1]{\underset{\mathcal{F}_n\!\!,#1}{\sim}}
\newcommand{\SV}{\underset{V}{\sim}}
\newcommand{\UF}[1]{\underset{\mathcal{F}\!\!,#1}{\prec}}
\newcommand{\ep}{\ensuremath{\varepsilon}}
\newcommand{\ho}{homomorphism}
\newcommand{\hos}{homomorphisms}
\newcommand{\fset}{\mathcal{F}}
\newcommand{\CCC}{\mathbb{C}}
\newcommand{\LL}{\mathcal{L}}
\newcommand{\LLL}{\mathcal{L}}
\newcommand{\KKK}{\mathcal{K}}
\newcommand{\EEE}{\mathcal{E}}
\newcommand{\HH}{\mathcal{H}}
\newcommand{\ot}{\otimes}
\newcommand{\fdr}{\mathrm{fdr}(A)}
\newcommand{\bN}{\bar{\mathbb{N}}}
\newcommand{\kkz}{KK(A,B)/\,\bar{0}}
\newcommand{\uk}{\underline{K}}
\newcommand{\homl}{\mathrm{Hom}_{\Lambda}(\uk(A),\uk(B))}
\newcommand{\homlp}{\mathrm{Hom}_{\Lambda}(\uk(A'),\uk(B))}
\newenvironment{rolist}{\begin{list}{(\roman{rocount})}
{\usecounter{rocount}\setlength{\itemsep}{0.0cm}\setlength{\parsep}
{0.0cm}}}{\end{list}} \setlength{\evensidemargin}{0in}
\setlength{\oddsidemargin}{0in} \setlength{\textwidth}{6in}
\newtheorem{theorem}{Theorem}[section]
\newtheorem{proposition}[theorem]{Proposition}
\newtheorem{lemma}[theorem]{Lemma}
\newtheorem{conjecture}[theorem]{Conjecture}
\newtheorem{corollary}[theorem]{Corollary}
\theoremstyle{remark}
\newtheorem{remark}[theorem]{Remark}
\newtheorem{examples}[theorem]{Examples}
\newtheorem{definition}[theorem]{Definition}
\newcommand{\Proof}{\medskip\noindent{\bf Proof.\ }}

\end{frontmatter}
\newcommand{\RG}{\ensuremath{R(\chi)}}
\newcommand{\MKB}{M(\mathcal{K}\otimes B)}
\newcommand{\bn}{\bar{\mathbb{N}}}
\newcommand{\SF}[1]{\underset{\mathcal{F}\!\!,#1}{\sim}}
\newcommand{\SFn}[1]{\underset{\mathcal{F}_n\!\!,#1}{\sim}}
\newcommand{\SV}{\underset{V}{\sim}}
\newcommand{\UF}[1]{\underset{\mathcal{F}\!\!,#1}{\prec}}
\newcommand{\ep}{\ensuremath{\varepsilon}}
\newcommand{\ho}{homomorphism}
\newcommand{\hos}{homomorphisms}
\newcommand{\fset}{\mathcal{F}}
\newcommand{\CCC}{\mathbb{C}}
\newcommand{\LL}{\mathcal{L}}
\newcommand{\LLL}{\mathcal{L}}
\newcommand{\KKK}{\mathcal{K}}
\newcommand{\EEE}{\mathcal{E}}
\newcommand{\HH}{\mathcal{H}}
\newcommand{\ot}{\otimes}
\newcommand{\fdr}{\mathrm{fdr}(A)}
\newcommand{\bN}{\bar{\mathbb{N}}}
\newcommand{\kkz}{KK(A,B)/\,\bar{0}}
\newcommand{\uk}{\underline{K}}
\newcommand{\homl}{\mathrm{Hom}_{\Lambda}(\uk(A),\uk(B))}
\newcommand{\homlp}{\mathrm{Hom}_{\Lambda}(\uk(A'),\uk(B))}

\section{Introduction}
Two $*$-homomorphisms $\varphi,\psi: A \to B$ are unitarily
equivalent if $u\varphi u^*=\psi$ for some unitary $u \in B$. They
are approximately unitarily equivalent, written
$\varphi\approx_{u}\psi$, if there is a sequence $(u_n)_{n \in
\mathbb{N}} $ of unitaries in $B$ such that $$\lim_{n \to \infty}
\|u_n \varphi (a) u_n^*- \psi (a)\|=0$$ for all $a \in A$.
Stable approximate unitary equivalence is a more elaborated concept
introduced in Def.~\ref{saue}.
 According to Glimm's
theorem,
 any non \mbox{type I} separable C*-algebra
 has uncountably many non unitarily equivalent
irreducible representations with the same kernel. In contrast, by
 Voiculescu's theorem, two
irreducible representations of a separable C*-algebra have the
same kernel if and only if they are approximately unitarily
equivalent. A comparison of the above results suggests that the
notion of
 unitary equivalence
is sometimes too rigid and that for certain purposes one can do
more things by  working with approximate unitary equivalence.
 This point of view is illustrated by Elliott's intertwining argument:
 if $\varphi:A \to B$ and $\psi:B \to A$ are unital $*$-homomorphisms
 between separable C*-algebras
 such that $\varphi\psi\approx_{u} id_B$ and $\psi\varphi\approx_{u} id_A$,
 then $A$ is isomorphic to $B$.
It is therefore very natural to study approximate unitary
equivalence of $*$-homomorphisms in a general context.

Two approximately unitarily equivalent $*$-homomorphisms
$\varphi,\psi: A \to B$ induce the same map on K-theory
 with coefficients, but they may have
different KK-theory classes. In order to handle this situation, R{\o}rdam
introduced the group $KL(A,B)$ as the quotient of $
\mathrm{Ext}(SA,B)^{-1}\cong KK(A,B)$ by the subgroup
$\mathrm{PExt}(K_{*-1}(A),K_*(B))$ of
 $\mathrm{Ext}(K_{*-1}(A),K_*(B))$
generated by pure  group extensions \cite{Ror:infsimple}. This required the
assumption that $A$ satisfies the universal coefficient theorem (UCT) of
\cite{RosSho:UCT}. Using a mapping cylinder construction, R{\o}rdam showed that
two approximately unitarily equivalent $*$-homomorphisms have the same class in
$KL(A,B)$. On the other hand, a  topology on  the Ext-theory groups was
considered by Brown-Douglas-Fillmore \cite{BroDouFil:Ext}, and shown to have
interesting applications in \cite{lgb:uct-qd} and \cite{Sa:relqd}. This
topology, called hereafter the Brown-Salinas topology, is defined via
approximate unitary equivalence of extensions. It was further investigated by
Schochet in
 \cite{Sc:fine1,Sc:fine2}  and by the author in \cite{Dad:ext}.
Schochet showed
 that the Kasparov product is continuous with respect to
the Brown-Salinas topology for K-nuclear separable C*-algebras. An important
idea from \cite{Sc:fine1,Sc:fine2} is that one can use the continuity of the
Kasparov product in order to transfer structural properties between
KK-equivalent C*-algebras. As it turns out, the subgroup
$\mathrm{PExt}(K_{*-1}(A),K_*(B))$ of
 $\mathrm{Ext}(SA,B)^{-1}$
 coincides with the closure of zero in the Brown-Salinas topology
 under the assumption that $A$ is nuclear and satisfies the UCT.
It is then quite natural  to define $KL(A,B)$ for arbitrary separable
C*-algebras as $\mathrm{Ext}(SA,B)^{-1}/\bar{0}$ as proposed by H. Lin in
\cite{Lin:auct}. Nevertheless, the study of $*$-homomorphisms
 from $A$ to $B$ via their class in $\mathrm{Ext}(SA,B)^{-1}$ is not
 optimal and leads to rather involved arguments as those in
 \cite{Lin:jot,Lin:auct} and \cite{Dad:ext} where
the Brown-Salinas topology of $\mathrm{Ext}(SA,B)^{-1}$ is related, in the
nuclear case, to stable approximate unitary equivalence of $*$-homomorphisms
from $A$ to $B$.

 Kasparov's
KK-theory admits several equivalent descriptions. This  deep
feature  enables one
 to choose  working with the picture
 that is most effective in a given situation.
  Similarly, there are several (and as we are going to see, equivalent)
  ways to
  introduce a topology on the
  KK-groups.  The Brown-Salinas topology was already mentioned.
  In a
recent important paper \cite{Pim:KK-top}, Pimsner defines a topology on the
equivariant graded KK-theory and proves the continuity of the Kasparov product
in full generality. The convergence of sequences in Pimsner's topology admits a
particularly nice and simple algebraic description
 which leads to major simplifications of the theory (see Lemma~\ref{Pimsner}).
However,  the previous  descriptions of the topology of $KK(A,B)$
 do not appear to be well adapted for the study of approximate unitary
 equivalence of $*$-\hos.

  In this paper we introduce  a
   topology on $KK(A,B)$ in terms of Cuntz pairs and
  approximate unitary equivalence.
   We then show that this topology
   coincides
  with Pimsner's topology
  (Thms.~\ref{mainconv} and ~\ref{2nd}).
  Our arguments rely  on a result of Thomsen \cite{Tho:onabex}
  and on our joint work with Eilers \cite{DadEil:AKK}.
Two $*$-homomorphisms
from $A$ to $B$ is the simplest instance of a Cuntz pair.
However, since in general the Kasparov group $KK(A,B)$
is not generated by $*$-homomorphisms from $A$ to $B$,
 it becomes necessary to work
with Cuntz pairs.
We revisit R{\o}rdam's group $KK(A,B)/\bar{0}$ in our general
setting and show
that it is a polish group (cf. \cite{Sc:fine1}) when endowed with the (quotient
of) Pimsner's topology for arbitrary separable C*-algebras (see
Prop.~\ref{23nd}). Along the way we show that the Brown-Salinas topology
coincides with Pimsner's topology (Cor.~\ref{P=BS}) and we give a series of
applications
 which include:

(i) two $*$-homomorphisms
 are stably approximately unitarily equivalent
 if and only if their KK-theory
 classes are equal modulo the closure of zero
 (see Cor.~\ref{for-homs},~\ref{for-ghoms}.)

 (ii) If  a separable C*-algebra $A$ satisfies the universal coefficient theorem in
KK-theory (UCT), then  $KK(A,B)/\,\bar{0}$ is homeomorphic to
$\mathrm{Hom}_{\Lambda}(\uk(A),\uk(B))$, where the latter group is endowed with
the topology of pointwise convergence (see Thm.~\ref{ppp}).
 Thus, in order to check that two KK-elements are close
 to each other, it suffices to verify that the maps
 they induce on the total K-theory group $\uk
(A)=\oplus_{n=0}^\infty K_*(A;\mathbb{Z}/n)$ agree on a sufficiently large
finite
 subset.

(iii) If a separable  exact residually finite dimensional C*-algebra satisfies
the UCT then it embeds in the UHF algebra of type $2^\infty$; see
Thm.~\ref{RRFFDD}. In particular the C*-algebra of a second countable amenable
locally compact maximally almost periodic group embeds in the UHF algebra of
type $2^\infty$.

(iv) We give a short proof of a theorem of H. Lin,
\cite{Lin:auct}, stating that two unital $*$-homomorphisms between
 Kirchberg C*-algebras are approximately unitarily
  equivalent if and only if their KL-classes coincide.
  This is used to show that a
   separable  nuclear C*-algebra satisfies the approximate
universal coefficient theorem of \cite{Lin:auct} if and only if it
satisfies the UCT
  (Thm.~\ref{sol}), answering a question
of H. Lin from \cite{Lin:auct}.

 For $A$ in the bootstrap category of \cite{RosSho:UCT}, one can  derive $(ii)$
from \cite{Sc:fine2} and \cite{DadLor:duke}. Its generalization to the
nonnuclear case is necessary in view of applications such as (iii). The latter
result was given a more complicated proof in an earlier preprint
\cite{Dad:emb-rfd}
 which is now superseded by the present paper.
A definition of the topology of $KK_{nuc}(A,B)$
 has also appeared there, but it became a  more useful tool after
 the emergence of \cite{Pim:KK-top}.
 The author is  grateful to M. Pimsner for
providing him with a draft of \cite{Pim:KK-top}.

\section{Metric structure}\label{metric section}
In this section we define an invariant pseudometric $d$ on $KK(A,B)$ which
makes $KK(A,B)$ a complete separable topological group. This is done by using a
description of $KK(A,B)$ based on Cuntz pairs and the asymptotic unitary
equivalence of \cite{DadEil:AKK}.

 The C*-algebras in this paper,
denoted by $A$, $B$, $C$,... will be assumed to be separable. We
only consider Hilbert $B$-bimodules, $E$, $F$,... that are
countably generated. The notation $H_B$ is reserved for the
canonical Hilbert $B$-bimodule obtained as the completion of
$\ell^2(\mathbb{N})\otimes_{alg} B$. As in \cite{Kas:cp} we
identify $M(B \otimes \KKK)$ with $\LLL(H_B)$. A unital
$*$-homomorphism $\gamma:A \to \LLL(H_B)$ is called \emph{unitally
absorbing} (for the pair of C*-algebras $(A,B)$) if for any unital
$*$-homomorphism $\varphi:A \to \LLL(H_B)$ there is a sequence of
unitaries $u_n \in \LLL(H_B,H_B \oplus H_B)$ such that for all $a
\in A$:

   (i) $\lim_{n \to \infty} \|u_n^* \left( \varphi(a)\oplus \gamma(a)
\right)u_n- \gamma(a)\|=0$

    (ii) $u_n^* \left( \varphi(a)\oplus \gamma(a) \right)u_n- \gamma(a)\in
\KKK(H_B)$

A $*$-homomorphism $\gamma:A \to \LLL(H_B)$ is called
\emph{absorbing} if its unitalization $\tilde{\gamma}:\tilde{A}
\to \LLL(H_B)$ is unitally absorbing.
 The theorems of
Voiculescu \cite{Voi:Weyl-vn} and Kasparov \cite{Kas:cp} exhibit
large classes of absorbing $*$-\hos. Thomsen \cite{Tho:onabex}
proved the existence of absorbing $*$-\hos\ for arbitrary
separable C*-algebras.

Let $\EEE_c(A,B)$ denote the set of all Cuntz pairs $(\varphi,\psi)$. They
consists of $*$-\hos\ $\varphi,\psi:A \to \LLL(H_B)$ such that
$\varphi(a)-\psi(a)\in \KKK(H_B)$ for all $a\in A$. It is was shown by Cuntz
that $KK(A,B)$ can be defined as the group of homotopy classes of Cuntz pairs.
In our joint work with Eilers we proved that $KK(A,B)$ can be realized  in
terms of proper asymptotic unitary equivalence classes of Cuntz pairs:
\begin{theorem}[\cite{DadEil:AKK}]\label{vanish}
Let $A$, $B$ be separable C*-algebras and let $(\varphi,\psi)\in \EEE_c(A,B)$
be a Cuntz pair. The following are equivalent:

(i)  $[\varphi,\psi]=0$  in $KK(A,B)$.

(ii) There is a $*$-homomorphism $\gamma:A \to \mathcal{L}(H_B)$
and there is a continuous unitary valued map
 $t \mapsto u_t \in 1+\mathcal{K}(H_B \oplus H_B)$,
 $t \in [0,\infty)$,
such that for all $a \in A$
\begin{equation}\label{conv}
\lim_{t \to \infty} \|u_t \left( \varphi (a)\oplus \gamma(a)
\right)u_t^*- \psi (a)\oplus \gamma(a)\|=0
\end{equation}

(iii) For any absorbing $*$-homomorphism $\gamma:A \to
\mathcal{L}(H_B)$ there is a continuous unitary valued map
 $t \mapsto u_t \in I+\mathcal{K}(H_B \oplus H_B)$,
 $t \in [0,\infty)$ satisfying \eqref{conv} for all $a \in A$.
\end{theorem}

This theorem suggests the following construction of a pseudometric
on $KK(A,B)$.

 Let $(a_i)_{i=1}^\infty$ be a dense sequence in the
unit ball of $A$. If $\varphi,\psi:A \to \LL(E)$ are
$*$-homomorphisms, we define
$$\delta_0(\varphi,\psi)=\sum_{i=1}^\infty
\frac{1}{2^i}\|\varphi(a_i)-\psi(a_i)\|,\,\, \text{and}$$
$$\delta_{\gamma}(\varphi,\psi)=
\inf\{\delta_0(\varphi\oplus \gamma,u(\psi \oplus \gamma)u^*): u \in
1+\KKK(E\oplus F)\,\text{unitary}\,\},$$ where $\gamma:A \to \LL(F)$ is an
absorbing $*$-\ho. One verifies immediately that
$\delta_\gamma(\varphi,\varphi)=0$,
$\delta_\gamma(\varphi,\psi)=\delta_\gamma(\psi,\varphi)$ and
$\delta_\gamma(\varphi,\eta)\leq\delta_\gamma(\varphi,\psi)+
\delta_\gamma(\psi,\eta).$ Moreover, if $\|\varphi_n(a)-\varphi(a)\|\to 0$ for
all $a \in A$, then $\delta_\gamma(\varphi_n,\psi)\to
\delta_\gamma(\varphi,\psi).$ If $\gamma_i:A \to \LL(F_i)$, $i=1,2$ are
$*$-\hos, then we write $\gamma_1\sim\gamma_2$ if there is a sequence of
unitaries $w_n \in \LL(F_1,F_2)$ such that for all $a \in A$
\begin{equation}\label{V}
\lim_{n\to \infty} \| w_n\gamma_1(a)w_n^*-\gamma_2(a)\|=0.
\end{equation}
\begin{lemma}\label{alpha}
If $\gamma_1\sim\gamma_2$, then
$\delta_{\gamma_1}(\varphi,\psi)=\delta_{\gamma_2}(\varphi,\psi).$
\end{lemma}
\begin{pf} If $w \in \LLL(F_1,F_2)$ is a
unitary, then $\delta_{\gamma_1}(\varphi,\psi)=\delta_{w\gamma_1
w^*}(\varphi,\psi)$, since conjugation by $1\oplus w$ maps
$1+\KKK(E\oplus F_1)$ onto $1+\KKK(E\oplus F_2)$. Thus
$\delta_{\gamma_1}(\varphi,\psi)=\delta_{w_n\gamma_1
w_n^*}(\varphi,\psi)\to\delta_{\gamma_2}(\varphi,\psi)$.
 \qed\end{pf} The assumption of Lemma~\ref{alpha} is automatically satisfied whenever $\gamma_i$ are absorbing $*$-\hos. Therefore
we can define $\delta(\varphi,\psi)=\delta_{\gamma}(\varphi,\psi)$ for some
absorbing $*$-\ho\ $\gamma$ and this definition does not depend on $\gamma$.

\begin{lemma}\label{beta}With notation as above

(a) If $w \in \LL(E,F)$ is a unitary, then $\delta(w\varphi
w^*,w\psi w^*)=\delta(\varphi,\psi),$

(b) If $\eta:A \to \LL(F)$ is a $*$-homomorphism, then
$\delta(\varphi,\psi)=\delta(\varphi \oplus \eta,\psi \oplus
\eta)=\delta(\eta \oplus \varphi,\eta \oplus\psi)$.
\end{lemma}
\begin{pf} For part (a) one argues as in the proof of the
previous lemma. For part (b) one uses the observation that
$\gamma\oplus \eta$ is absorbing whenever $\gamma$ is absorbing and
part (a). \qed\end{pf} If $\varphi,\psi:A \to \LLL(E)$ are $*$-\hos,
we write $(\varphi)\approx (\varphi')$ if there is   a sequence of
unitaries $u_n \in 1+\KKK(E)$ such that $\lim_{n\to
\infty}\|u_n\varphi(a)u_n^*-\psi(a)\|=0$ for all $a \in A$.

\begin{lemma}\label{gamma}
Let $\varphi,\psi:A \to \LLL(E)$ and $\varphi',\psi':A \to
\LLL(E)$ be $*$-\hos. Assume that $(\varphi)\approx (\varphi')$
and $(\psi)\approx (\psi')$. Then
$\delta(\varphi,\psi)=\delta(\varphi',\psi')$.
\end{lemma}
\begin{pf} This is an immediate consequence of the definition of
$\delta$ and the observation that if $u\in 1+\KKK(E)$ is a unitary,
then $\delta(u\varphi u^*,\psi)=\delta(\varphi,\psi)$. \qed\end{pf}

 We are now ready to introduce a pseudometric $d$ on
$\EEE_c(A,B)$. A pseudometric  satisfies all the properties of a
metric except that $d(x,y)=0$ may not imply $x=y$.
\begin{definition}
$d((\varphi,\psi),(\varphi',\psi'))=\delta(\varphi \oplus \psi',\psi
\oplus \varphi').$
\end{definition}
\begin{lemma}\label{crucial}
If $x,x'\in \EEE_c(A,B)$ and $[x]=[x']$ in $KK(A,B)$ then
$d(x,x')=0$.
\end{lemma}
\begin{pf}If $x=(\varphi,\psi)$ and $x=(\varphi',\psi')$ then
$[x]-[x']=[\varphi\oplus\psi',\psi\oplus\varphi']=0$. By Theorem
\ref{vanish}  this implies
$\delta(\varphi\oplus\psi',\psi\oplus\varphi')=0$ hence
$d(x,x')=0$. \qed\end{pf}
\begin{proposition}\label{delta}
$d$ is a pseudometric on $\EEE_c(A,B)$ that descends to an invariant
pseudometric on $KK(A,B)$ (denoted again by $d$!).
\end{proposition}
\begin{pf}
First we show that $d$ is a pseudometric on $\EEE_c(A,B)$. Let
$x=(\varphi,\psi), x'=(\varphi',\psi')\in \EEE_c(A,B)$. Then
$d(x,x)=0$ by Lemma~\ref{crucial}.  The equality $d(x,x')=d(x',x)$
is equivalent to $\delta(\varphi \oplus \psi',\psi \oplus
\varphi')=\delta(\varphi' \oplus \psi ,\psi' \oplus \varphi)$. The
latter equality follows from Lemma~\ref{beta}(a) with $w$ a
permutation unitary and the symmetry of $\delta$. In order to verify
the triangle inequality for $d$, we first recall that if
$\alpha,\alpha',\alpha'':A \to \LLL(E)$ then
\begin{equation}\label{tria}
\delta(\alpha,\alpha')+\delta( \alpha',\alpha'')\geq
\delta(\alpha,\alpha'').
\end{equation}
 Let
$x''=(\varphi'',\psi'')\in \EEE_c(A,B)$. The inequality
$d(x,x')+d(x',x'')\geq d(x,x'')$ is equivalent to
\begin{equation}\label{target}
    \delta(\varphi\oplus\psi',\psi\oplus\varphi')+
    \delta(\varphi'\oplus\psi'',\psi'\oplus\varphi'') \geq
    \delta(\varphi\oplus\psi'',\psi\oplus\varphi'')
\end{equation}
By Lemma~\ref{beta}
$$\delta(\varphi\oplus\psi'',\psi\oplus\varphi'')=
\delta(\varphi\oplus\psi''\oplus\psi',\psi\oplus\varphi''\oplus\psi')=
\delta(\varphi\oplus\psi'\oplus\psi'',\psi\oplus\psi'\oplus\varphi'')$$
and the latter term less than or equal to
$\delta(\varphi\oplus\psi'\oplus\psi'',\psi\oplus\varphi'\oplus\psi'')
+\delta(\psi\oplus\varphi'\oplus\psi'',\psi\oplus\psi'\oplus\varphi'')$
by \eqref{tria}. Finally
$\delta(\varphi\oplus\psi'\oplus\psi'',\psi\oplus\varphi'\oplus\psi'')=
\delta(\varphi\oplus\psi',\psi\oplus\varphi')$ and
$\delta(\psi\oplus\varphi'\oplus\psi'',\psi\oplus\psi'\oplus\varphi'')=
\delta(\varphi'\oplus\psi'',\psi'\oplus\varphi'')$ by
Lemma~\ref{beta}. This proves the inequality \eqref{target}.

Next we are going to verify that $d$ descends to a metric on $KK(A,B)$. By
symmetry, it suffices to prove that if $x,x',x''\in \EEE_c(A,B)$ and
$[x']=[x'']$ in $KK(A,B)$, then $d(x,x'')\leq d(x,x')$. By Lemma~\ref{crucial},
$d(x',x'')=0$. Since $d$ is a pseudometric, $d(x,x'')\leq
d(x,x')+d(x',x'')=d(x,x')$.

It remains to verify the invariance of the pseudometric. We show
that $d(x\oplus y,x'\oplus y)=d(x,x')$ for all $x,x',y \in
\EEE_c(A,B)$. Let $\hat{d}([x],[x'])=d(x,x')$ denote (temporarily)
the induced metric on $KK(A,B)$. We claim that
$d(x,x')=\hat{d}([x]-[x'],0)$, which implies the invariance of $d$.
To verify the claim note that if $x=(\varphi,\psi)$ and
$x=(\varphi',\psi')$ then
$d(x,x')=\delta(\varphi\oplus\psi',\psi\oplus\varphi')$ by
definition, and
$\hat{d}([x]-[x'],0)=d((\varphi\oplus\psi',\psi\oplus\varphi'),(0,0))=
\delta(\varphi\oplus\psi',\psi\oplus\varphi')$. \qed\end{pf}
\begin{proposition}\label{23nd} Let $A$ be $B$ be separable C*-algebras.
The topology of $KK(A,B)$ defined by the pseudometric $d$ satisfies the
 second axiom of countability.
 If $\bar{0}$ denotes the closure of zero, then $\kkz$ is a polish group.
\end{proposition}
\begin{pf}  By a result of Thomsen \cite[Thm.~3.2]{Tho:onabex}, every
element of $KK(A,B)$ is represented by a Cuntz pair $(\alpha,\gamma)$, where
$\gamma:A \to \LLL(H_B)$ is a fixed  absorbing $*$-\ho. Therefore the image of
each map $\alpha$ is contained in the separable C*-algebra
$\gamma(A)+\KKK(H_B)$. This shows that
  the topology of $KK(A,B)$ satisfies the second axiom of
countability.

 Next we prove the completeness of $KK(A,B)$. Let
$(x_n)$ be a Cauchy sequence in $\EEE_c(A,B)$ where
$x_n=(\alpha_n,\gamma)$ with $\gamma:A \to \LLL(H_B)$ as above. This
means that
$d(x_n,x_m)=\delta(\alpha_n\oplus\gamma,\gamma\oplus\alpha_m)\to 0$
as $m,n \to \infty$. Since
$\delta(\alpha_m\oplus\gamma,\gamma\oplus\alpha_m)=d(x_m,x_m)=0$,
 we have
$\delta(\alpha_n\oplus\gamma,\alpha_m\oplus\gamma)\to 0$ as $m,n \to \infty$.
Since $[\alpha_n,\gamma]=[\alpha_n\oplus\gamma,\gamma\oplus\gamma]$ in
$KK(A,B)$, after replacing $\alpha_n$ by $\alpha_n\oplus\gamma$, we may assume
that $\delta(\alpha_n,\alpha_m)\to 0$ as $m,n \to \infty$. After passing to a
subsequence of $(\alpha_n)$, if necessary, we find a sequence of unitaries
$u_n\in 1+\KKK(H_B)$ such that
$\delta_0(\alpha_n,u_{n+1}\alpha_{n+1}u_{n+1}^*)<1/2^n$. Define
$\alpha'_n(a)=(u_2\cdots u_n)\alpha_n(a) (u_2\cdots u_n)^*$ and note that
$(\alpha'_n)$ is a Cauchy sequence in $\mathrm{Hom}(A,\LLL(H_B))$ since
$\delta_0(\alpha'_n,\alpha'_{n+1})<1/2^n$. Since $\mathrm{Hom}(A,\LLL(H_B))$ is
complete, $(\alpha'_n)$ converges to a $*$-\ho\ $\alpha$ with the property that
$\alpha(a)-\gamma(a) \in \KKK(H_B)$ since $\alpha'_n(a)-\gamma(a) \in
\KKK(H_B)$ for all $a \in A$. It follows that
$[\alpha_n,\gamma]=[\alpha'_n,\gamma]$ converges to $[\alpha,\gamma]$ in
$KK(A,B)$. \qed\end{pf} Proposition~\ref{23nd} does not follow from
\cite{Sc:fine1} since we do not assume $A$ to be K-nuclear and we are working a
priori with a different topology.
\section{Approximate unitary equivalence and the topology of $KK(A,B)$}\label{two}
In this section we show that the approximate unitary equivalence
of Cuntz pairs can be expressed in KK-theoretical terms, see
Theorem~\ref{mainconv}. Consequently, the topology of $KK(A,B)$
defined by $d$ coincides with  Pimsner's topology, see
Theorem~\ref{2nd}.
  In the final part we apply these results to
$*$-\hos.

Let $\bn=\{1,2,\dots\}\cup\{\infty\}$ denote the one-point compactification of
the natural numbers. We say that a topology on the $KK$-theory groups satisfies
Pimsner's condition if the convergence of sequences  is characterized as
follows. A sequence $(x_n)$ in $KK(A,B)$ converges to $x_\infty$ if and only if
there is $y \in KK(A,C(\bar{\mathbb{N}})\otimes B)$ with $y(n)= x_n$ for $n \in
\mathbb{N}$ and $y(\infty)=x_\infty$. Clearly  a topology which satisfies the
first axiom of countability and Pimsner's condition is unique. Pimsner made the
following crucial observation.

\begin{lemma}\cite{Pim:KK-top}\label{Pimsner}
If a topology on the $KK$-groups satisfies the first axiom of countability and
Pimsner's condition, then the Kasparov product is jointly continuous with
respect to that topology.
\end{lemma}
\begin{pf} By the functoriality of the cup product of Kasparov
\cite[2.14]{Kas:inv}, if $y \in KK(A,C(\bar{\mathbb{N}})\otimes B)$ and $z \in
KK(B,C(\bar{\mathbb{N}})\otimes C)$, then the image $w \in
KK(A,C(\bar{\mathbb{N}})\otimes C)$  of the cup product $y \otimes_B z \in
KK(A,C(\bar{\mathbb{N}}\times \bar{\mathbb{N}})\otimes C)$ under the diagonal
map satisfies $w(n)=z(n)\otimes y(n)$ for all $n \in \bN$.\qed\end{pf}

 We need some notation. Let $\fset \subset A$ be a finite subset
and let $\ep >0$. If $\varphi:A \to \LLL_{B }(E)$ and $\psi:A \to \LLL_{B }(F)$
are two contractive completely positive maps, we write $\varphi \UF{\ep} \psi$
if there is an isometry $v\in \LLL_{B}(E, F) $ such that
$\|\varphi(a)-v^*\psi(a)v\|< \ep$ for all $a \in \fset$. If $v$ can be taken to
be a unitary then we write  $\varphi \SF{\ep} \psi$. We write $\varphi \prec
\psi$ (respectively $\varphi \sim \psi$) if
 $\varphi \UF{\ep} \psi$ (respectively $\varphi \SF{\ep} \psi$)
 for all finite sets $\fset$ and $\ep>0$. Note that if $\varphi \UF{\ep_1} \psi$
 and $\psi \UF{\ep_2} \gamma$, then $\varphi \UF{\ep_1+\ep_2} \gamma$.

\begin{proposition}\label{abs-nuc}
Let $A$, $B$, $C$ be separable C*-algebras such that $B$ stable and $C$ is
unital and nuclear. If $\gamma:A \to M(B)$ is an absorbing $*$-homomorphism for
$(A,B)$, then $\Gamma:A \to M(B\otimes C)$, $\Gamma(a)=\gamma(a)\otimes 1_C$,
is an absorbing $*$-homomorphism for $(A,B\otimes C)$.
\end{proposition}

\begin{pf} If $B=\KKK$ this is essentially Kasparov's absorption
theorem \cite{Kas:cp}. By \cite[Thm.~2.13]{DadEil:class} it suffices to prove
that for any finite subset $\fset \subset A$, any $\ep
>0$ and any completely positive contraction  $\sigma:A \to B\ot C$
we have $\sigma \UF{\ep} \Gamma$. Since $\gamma$ is an absorbing
$*$-homomorphism for $(A,B)$, we have $\Phi \prec \gamma$ and hence $\Phi \ot
1_C \prec \gamma \ot 1_C=\Gamma$ for any completely positive contraction
$\Phi:A \to B$. Therefore it is enough to show that $\sigma \UF{\ep} \Phi \ot
1_C$ for some completely positive contraction  $\Phi:A \to B$. Since $C$ is
nuclear, as a consequence of Kasparov's theorem, $id_C \prec \delta \ot 1_C$
where $\delta:C \to \LLL(H)$ is a unital faithful representation with
$\delta(C)\cap \KKK(H)=\{0\}$. Therefore there is sequence of isometries
$v_n\in \LL_C(C,H_C)$ with
$$\lim_{n \to \infty} \|c-v_n^*(\delta(c)\ot 1_C)v_n\|=0$$ for all
$c \in C$. Since $H_C$ is the closure of $\oplus_{n=1}^\infty C$ one can
perturb each $v_n$ to a $C$-linear isometry   $v_n :C \to C^{k(n)}\subset H_C$.
Therefore if $\delta_n:C \to M_{k(n)}(\CCC)$ denotes the completely positive
contraction obtained by compressing $\delta$ to the subspace
$\mathbb{C}^{k(n)}$ of $H$, we have
$$\lim_{n \to \infty} \|c-v_n^*(\delta_n(c)\ot
1_C)v_n\|=0$$ for all $c \in C$. If we set $V_n=\mathrm{id}_B \ot
v_n \in \LL_{B\ot C} (B\ot C, (B\ot C)^{k(n)})$ and
$\Delta_n=\mathrm{id}_B\ot \delta_n : B\ot C \to B\ot M_{k(n)}(
\CCC)$, then
$$\lim_{n \to \infty} \|x-V_n^*(\Delta_n(x)\ot
1_C)V_n\|=0$$ for all $x \in B \ot C$. Consequently
\begin{equation}\label{q}
\lim_{n \to \infty} \|\sigma(a)-V_n^*(\Delta_n(\sigma(a))\ot 1_C)V_n\|=0
\end{equation} for all $a \in A$. Note that $\Phi_n=\Delta_n \sigma:A \to
M_{k(n)}(B)\cong B$ is a completely positive contraction. From \eqref{q} we see
that $\sigma \UF{\ep} \Phi_n \ot 1_C$ for some large enough $n$ and this
concludes the proof. \qed\end{pf}

\begin{theorem}\label{mainconv}
Let $A$, $B$ be separable C*-algebras and let $(\varphi_n,\psi_n)_{n\in
\mathbb{N}}$  be a sequence of Cuntz pairs in $ \mathcal{E}_c(A,B)$. The
following are equivalent:

(i) There is $y \in KK(A,C(\bar{\mathbb{N}})\otimes B)$ such that $y(n)=
[\varphi_n,\psi_n]$ for $n \in \mathbb{N}$ and $y(\infty)=0$.

(ii) For any absorbing $*$-homomorphism $\gamma:A \to
\mathcal{L}(H_B)$ there is a sequence of unitaries $u_n \in
1+\mathcal{K}(H_B \oplus H_B)$ such that for all $a \in A$
\begin{equation}\label{converge}
\lim_{n \to \infty} \|u_n \left( \varphi_n(a)\oplus \gamma(a) \right)u_n^*-
\psi_n(a)\oplus \gamma(a)\|=0
\end{equation}
(iii) The sequence $[\varphi_n,\psi_n]$ converges to zero in
$(KK(A,B),d)$.
\end{theorem}
\begin{remark}  It is
easy to verify that condition (ii) is equivalent to asking  that
there is some $*$-homomorphism $\gamma:A \to \mathcal{L}(H_B)$ and
there is a sequence of unitaries $u_n \in I+\mathcal{K}(H_B \oplus
H_B)$
 satisfying \eqref{converge}
for all $a \in A$. This is very similar to the proof of (ii) $
\Leftrightarrow$ (iii) of Theorem~\ref{vanish}.\end{remark}
\begin{pf}
Given two sequence of $*$-homomorphisms $\varphi_n,\psi_n:A \to
\LL(E_n)$,
 we write   $(\varphi_n)_n\approx
(\psi_n)_n$ if there is a sequence of unitaries $u_n \in
1+\KKK(E_n)$ such that
$$\lim_{n\to \infty}\|u_n\varphi_n(a)u_n^*-\psi_n(a)\|=0$$
for all $a \in A$.  With this notation, the
  condition ~\eqref{converge} reads $(\varphi_n\oplus \gamma)_n\approx
  (\psi_n\oplus\gamma)_n$. It is easy to verify that $\approx$ is an
equivalence relation and that $(\varphi_n\oplus
\varphi'_n)_n\approx (\psi_n\oplus \psi'_n)_n$ whenever
$(\varphi_n)_n\approx (\psi_n)_n$ and
  $(\varphi'_n)_n\approx (\psi'_n)_n$.

  We identify
  $\LLL(H_B)$ with $M(\KKK \otimes B  )$ and $\KKK(H_B)$ with $\KKK \otimes B$.
   Therefore the set $\EEE_c(A,B)$ consists of
  pairs of $*$-\hos\ $(\varphi,\psi):A \to M(\KKK\ot B)$ such that
  $\varphi(a)-\psi(a) \in \KKK\ot B$ for all $a \in A$.
  Since $M(\KKK\ot B\ot C(\bN))\equiv C_{s}(\bN,M(\KKK\ot B))$
  (the set of strictly continuous functions from $\bn$ to $M(\KKK\ot B)$)
  and $\KKK\ot B\ot
  C(\bN)=C(\bN,\KKK\ot B)$,
   an element $(\Delta,\Gamma) \in \EEE_c(A,B\ot C(\bN))$
  is completely determined by a family $(\delta_n,\gamma_n)_{n \in \bN}\subset \EEE_c(A,B)$
  such that \begin{equation}\label{tobes}\lim_{n \to
  \infty}\delta_n(a)=\delta_\infty(a),\quad
\lim_{n \to
  \infty}\gamma_n(a)=\gamma_\infty(a)\end{equation}
  in the strict topology of $M(\KKK\otimes B)$
  and such that
  $$\lim_{n \to
  \infty}(\delta_n(a)-\gamma_n(a))=\delta_\infty(a)-\gamma_\infty(a)$$
   in the norm topology, for all $a \in A$.
   By \cite[Thm.~3.2]{Tho:onabex},  each element of $KK(A,B)$
   is represented by a pair $(\delta,\gamma)\in \EEE_c(A,B)$ where
   $\gamma:A \to M(\KKK\ot B)$ is any given  absorbing $*$-\ho. In view of
   Proposition~\ref{abs-nuc}, if $y \in KK(A,B \ot C(\bN))$, then we
   can write
   $y=[\Delta,\Gamma]$ where $\Gamma=\gamma \ot 1_{C(\bN)}$ and $\gamma:A \to M(\KKK\otimes B)$ is a
   fixed  absorbing $*$-\ho\ for $(A,B)$. In other words $\Gamma$ is
   given by a constant family $(\gamma_n)_{n \in \bn}$ with
   $\gamma_n=\gamma$.
   A crucial consequence of our choice of $\Gamma$ is that
   $\delta_n(a)-\delta_\infty(a)\in \KKK\ot B$ for all $a \in A$, since it is equal to
   $(\delta_n(a)-\gamma(a))-(\delta_\infty(a)-\gamma(a))$
   and therefore
   \begin{equation}\label{cccrucial}
   \lim_{n \to
  \infty}\|\delta_n(a)-\delta_\infty(a)\|=0
  \end{equation}
  for all $a \in A$. Therefore we are able to pass
  from strict convergence in \eqref{tobes} to norm convergence in \eqref{cccrucial}.
After this preliminary discussion we proceed with the proof of the
theorem. The equivalence (ii)$\Leftrightarrow$ (iii) follows
immediately from the definition of $d$ and the separability of $A$.

(i) $\Rightarrow$ (ii) It is convenient to consider first the
situation when $(\varphi_n,\psi_n)$ is a sequence of Cuntz pairs
where the second component $\psi_n$ is fixed for all $n$ and equal
to some  absorbing $*$-homomorphism  $\gamma$ as above. By
assumption there is $y \in KK(A,B \ot C(\bN))$ such that
$y(n)=[\varphi_n,\gamma]$ and $y(\infty)=0$. Write
$y=[\Delta,\Gamma]$  as above. Therefore
$[\delta_n,\gamma]=[\varphi_n,\gamma]$ hence
$[\delta_n,\varphi_n]=0$ and $[\delta_\infty,\gamma]=0$. Using
Theorem~\ref{vanish} we obtain
$$(\delta_n\oplus \gamma)_n\approx (\varphi_n\oplus \gamma)_n,\quad
(\delta_\infty\oplus \gamma)_n\approx (\gamma\oplus \gamma)_n.$$
In view of \eqref{cccrucial} this gives
\begin{equation}\label{prelim}
    (\varphi_n \oplus \gamma)_n\approx (\gamma \oplus \gamma)_n.
\end{equation}
We now proceed with the general case  with $(\varphi_n,\psi_n)$ as in (i).
Using \cite[Thm.~3.2]{Tho:onabex} again, we find a sequence $(\gamma_n,\gamma)
\in \EEE_c(A,B)$ with $[\gamma_n,\gamma]=[\varphi_n,\psi_n]$ and $\gamma$
absorbing. Since $[\varphi_n \oplus \gamma ,\psi_n \oplus \gamma_n]=0$, by
Theorem~\ref{vanish}  we obtain $(\varphi_n \oplus \gamma \oplus
\gamma)_n\approx (\psi_n \oplus \gamma_n \oplus \gamma)_n$.
 By the
first part of the proof, we have $(\gamma_n \oplus
\gamma)_n\approx (\gamma \oplus \gamma)_n$. Altogether this gives
$(\varphi_n \oplus \gamma \oplus \gamma)_n\approx (\psi_n \oplus
\gamma \oplus \gamma)_n$. Since $\gamma$ is  absorbing, $(\gamma
\oplus \gamma)_n\approx (\gamma)_n$ hence we obtain (ii):
$(\varphi_n \oplus \gamma)_n\approx (\psi_n \oplus \gamma)_n$.

(ii) $\Rightarrow$ (i) Replacing $\varphi_n$ by $\varphi_n \oplus \gamma$ and
$\psi_n$ by $\psi_n \oplus \gamma$ we may assume that there are unitaries $u_n
\in I+\KKK(H_B)$ such that  $\lim_{n \to
\infty}\|u_n\varphi_n(a)u_n^*-\psi_n(a)\|=0$ for all $a \in A$. Since
$(u_n\varphi_n u_n^*,\psi_n)$ and $( \varphi_n,\psi_n)$ have the same KK-class,
after replacing $\varphi_n$ by $u_n\varphi_n u_n^*$ we may further assume that
$\lim_{n \to \infty}\| \varphi_n(a)-\psi_n(a)\|=0$. Since both  $\varphi_n$ and
$\gamma$ are  absorbing, there is a sequence of unitaries $w_n \in \LL(H_B)$
such that $ w_n \varphi_n(a)w_n^*-\gamma(a)\in \KKK(H_B)$ and $\lim_{n \to
\infty}\|w_n \varphi_n(a)w_n^*-\gamma(a)\|=0$. Define $*$-\hos\ $\Phi,\Psi:A
\to M(\KKK\otimes B \otimes C(\bN))$ by setting $\Phi_n=w_n \varphi_n w_n^*$,
$\Phi_\infty=\gamma$, $\Psi_n=w_n \psi_n w_n^*$ and $\Psi_\infty=\gamma$. The
family $(\Phi,\Psi)=(\Phi_n,\Psi_n)_{n \in \bN}$ defines an element $y$ of
$KK(A,C(\bN)\ot B)$ such that $y(n)=[\Phi_n,\Psi_n]=[\varphi_n,\psi_n]$ for $n
\in \mathbb{N}$ and $y(\infty)= [\Phi_\infty,\Psi_\infty]=[\gamma,\gamma]=0$.
\qed\end{pf} We collect the previous results of the section in the following
form.
\begin{theorem}\label{2nd} Let $A$ be $B$ be separable C*-algebras.
The topology of $KK(A,B)$ defined by the pseudometric $d$ is
separable and complete. A sequence $(x_n)_{n=1}^\infty$ converges
to
  $x_\infty$
 in $KK(A,B)$ if and only if there is $y \in KK(A,C(\bar{\mathbb{N}})\otimes B)$
 with $y(n)= x_n$ for all $n \in \bn$. Therefore the topology defined by
 $d$ satisfies Pimsner's condition and hence the Kasparov product is continuous.
 The topology defined by
 $d$ coincides with Pimsner's topology.
\end{theorem}
\begin{pf} The first part follows from Proposition~\ref{23nd} and
Theorem~\ref{mainconv}. The second part follows from
Lemma~\ref{Pimsner}.
 \qed\end{pf}

Let us see how the previous results
 can be applied to $*$-homomorphisms. The definition of
stable approximate unitary equivalence for two $*$-\hos\
$\varphi,\psi :A \to B$
 is not quite straightforward. A naive definition that would require
 approximate unitary equivalence after taking direct sums with
$*$-\hos\ would not be satisfactory, due to a possible small supply of
$*$-homomorphisms from $A$ to $B$.

\begin{definition}\label{saue} Let $A$, $B$ be  separable C*-algebras.
Two $*$-\hos\ $\varphi,\psi :A \to B$ are called
stably approximately unitarily equivalent if there is
a sequence of unitaries
$v_n \in 1+\mathcal{K}(B \oplus H_B)$
 and an absorbing $*$-homomorphism $\gamma:A \to M(B \otimes \KKK)$
such that for all $a \in A$
\begin{equation}\label{saue-hom}
\lim_{n \to \infty} \|v_n \left( \varphi(a)\oplus \gamma(a) \right)v_n^*-
\psi(a)\oplus \gamma(a)\|=0
\end{equation}
\end{definition}
From Theorem~\ref{mainconv} we obtain:
\begin{corollary}\label{for-ghoms}
Let $A$, $B$ be  separable C*-algebras. Two $*$-\hos\ $\varphi,\psi :A \to B$
are stably approximately unitarily equivalent if and only if
$[\varphi]-[\psi]\in\bar{0}$ in $KK(A,B)$, if and only if
$d([\varphi],[\psi])=0.$
\end{corollary}

This result becomes more useful when there are many $*$-\hos\ from $A$ to $B$
or matrices over $B$. For illustration, we generalize \cite[Thm.~5.1]{Dad:ext}
and \cite[Thm.~3.9]{Lin:auct}.
 Let  $A$ and $B$ be unital
separable C*-algebras such that either $A$ or $B$ is nuclear. Assume
 that there is a sequence of unital $*$-homomorphisms
$\eta_n:A \to M_{k(n)}(B)$  such that for all nonzero $a\in A$ the closed
two-sided ideal of $B\otimes \KKK$ generated by $\{\eta_n(a):n \in
\mathbb{N}\}$ is equal to $B\otimes \KKK$. We will also assume that each
$\eta_n$ appears infinitely many times in the sequence $(\eta_n)$.

\begin{corollary}\label{for-homs}
Let  $A$ and $B$ be unital separable C*-algebras such that either $A$ or $B$ is
nuclear. Assume that $(\eta_n)$ is as above and let $\varphi$, $\psi$ be two
unital $*$-homomorphisms from $A$ to $B$. Then $[\varphi]-[\psi]\in\bar{0}$
 if and only if there exist a sequence on integers $(m(n))$ and
unitaries $(u_n)$ in matrices over $B$ (of appropriate size)
 such that
\begin{equation}\label{converges-hom}\lim_{n \to \infty}
\|u_n (\varphi (a)\oplus \gamma_n(a) )u_n^*- \psi (a)\oplus
\gamma_n(a)\|=0\end{equation} for all $a \in A$, where $\gamma_n=\eta_1\oplus
\dots \oplus \eta_{m(n)}$.
\end{corollary}
\begin{pf}
We verify only the nontrivial implication $(\Rightarrow)$. To simplify
notation, we give the proof in the case when all $k(n)=1$, i.e. $\eta_n:A \to
B$. The condition $[\varphi]-[\psi]\in\bar{0}$ is equivalent to the condition
that $0$ belongs to the closure of $[\varphi]-[\psi]$. If $\gamma:A \to M(B
\otimes \KKK)$ is defined  by
$$\gamma(a)=\mathrm{diag}(\eta_1(a),\eta_2(a),\cdots),$$
then $\gamma$ is a unitally absorbing representation by a result of
\cite{Ell-Ku}. If $\widehat \gamma:A \to \LLL(H_B\oplus H_B)\cong
\LLL(H_B)$ is defined by $\widehat \gamma=\gamma\oplus 0$, then
$\widehat \gamma$ is absorbing. By Theorem~\ref{mainconv} there is a
sequence of unitaries $v_n \in 1+\mathcal{K}(B \oplus H_B)$ such
that for all $a \in A$
\begin{equation}\label{converge-hom}
\lim_{n \to \infty} \|v_n \left( \varphi(a)\oplus \widehat\gamma(a)
\right)v_n^*- \psi(a)\oplus \widehat\gamma(a)\|=0
\end{equation}
If $e_m=1_B \oplus\cdots\oplus1_B$ ($m$-times), then
$\lim_{m \to \infty} \|[v_n ,e_m]\|=0$ for all $m$.
For each $n$ let $m(n)$ be such that $\|[v_n ,e_{m(n)}]\|<1/n$.
By functional calculus, there are unitaries
$u_n \in M_{m(n)}(B)$ with
$\lim_{n\to \infty} \|u_n -e_{m(n)}v_ne_{m(n)}\|=0$.
With these choices we derive \eqref{converges-hom} by
 compressing in \eqref{converge-hom} by $e_{m(n)}$.
\qed\end{pf} The following result is derived by a similar argument.
 \begin{corollary}\label{nice}
Let $A$, $B$ and $(\eta_n)$ be as in Cor.~\ref{for-homs} and let $\varphi$,
$\psi$ be two unital $*$-homomorphisms from $A$ to $B$. For any finite subset
$\fset$ of $A$ and any $\varepsilon>0$ there is $\delta>0$ such that if
$d([\varphi],[\psi])<\delta$ then $\varphi\oplus \gamma_n \SF{\ep}
\psi\oplus\gamma_n$ for some $n$.
\end{corollary}
\section{The UCT, K-theory with coefficients
 and applications}\label{UCT}
Let $A$ and $B$ be separable C*-algebras. The total K-theory group $\uk
(A)=\oplus_{n=0}^\infty K_*(A;\mathbb{Z}/n)$ has a natural action of the
Bockstein operations $\Lambda$ of \cite{Sc:bock}. In this section we show that
if $A$ satisfies the UCT, then $\kkz$ is isomorphic as a topological group with
$\mathrm{Hom}_{\Lambda}(\uk(A),\uk(B))$ endowed with the topology of
 pointwise convergence. This is extremely useful
 since in order to check that two KK-elements are close
 to each other it suffices to show that the maps
 they induce on $\uk(-)$ agree on some (sufficiently large) finite subset.
 By a result of J.L.
Tu \cite{Tu:BC}, the C*-algebra of an a-T-menable locally compact
second countable groupoid with Haar system satisfies the UCT. This
shows that there are large natural classes of non-nuclear
C*-algebras satisfying the UCT. As an application we show that the
C*-algebra of  a second countable amenable locally compact
maximally almost periodic group embeds in the UHF algebra of type
$2^\infty$.

  If $d_*$ is the
metric on $\uk(B)$ with $d_*(x,y)=1$ for $x\neq y$, then $\homl$
becomes a polish group with respect to the metric
$$d(\mu,\nu)=\sum_{n=1}^\infty
\frac{1}{2^n}d_*(\mu(x_n),\nu(x_n)),$$ where $\{x_1,x_2,\dots\}$
is an enumeration $\uk(A)$.

  A separable C*-algebra satisfies the UCT of \cite{RosSho:UCT}
  if and only if is
 KK-equivalent to a commutative C*-algebra, if and only it
 satisfies the following universal multi-coefficient exact
 sequence of \cite{DadLor:duke}:
 \begin{equation}\label{UMCT}
0\to \mathrm{PExt}(K_{*-1}(A),K_*(B)) \to KK(A,B)
\stackrel{\Gamma}{\rightarrow}
 \mathrm{Hom}_\Lambda(\uk (A),\uk (B) )\to 0.
 \end{equation}
 Here $\mathrm{PExt}$ stands for the subgroup of $\mathrm{Ext}$ corresponding
 to pure extensions. We refer the reader to the monograph \cite{sc:primer}
 for an excellent introduction to $\mathrm{PExt}$.
 The map $\Gamma$ is induced by the Kasparov product
 and therefore is continuous.
 This is also easily seen directly since if
two projections are close to each other then they have the same
K-theory class.

If $x\in KK(A,B)$ we denote $\Gamma(x)$ by $\underline{x}$. The
following result can be deduced from \cite{Sc:fine2} for nuclear
C*-algebras $A$ in the bootstrap category of \cite{Sc:II}, modulo
the identification of Pimsner's topology with the Brown-Salinas
topology. The idea of using the continuity of the Kasparov product
in its proof is borrowed from \cite{Sc:fine2}.

\begin{theorem}\label{ppp} Let $A$ and $B$ be separable C*-algebras and assume
that $A$ satisfies the UCT. Then

(a) $x_n \to x$ in $KK(A,B)$ if and only if $\underline{x_n}\to
\underline{x}$ in $\homl$.

(b) The map $\kkz \to \homl$ is an isomorphism of topological
groups. In particular $\kkz$ is totally disconnected.
\end{theorem}
\begin{pf}
Part (a) is an immediate consequence of (b).  Since the Kasparov
product is continuous, multiplication by a KK-invertible element
$y\in KK(A,A')$ induces a commutative diagram
$$ \xymatrix{
{KK(A',B)}\ar[r]\ar[d]&
{\homlp}\ar[d]\\
{KK(A,B)}\ar[r]&{\homl} }$$ where the horizontal maps are
continuous and the vertical maps are homeomorphisms. Therefore,
 after replacing $A$ by a
KK-equivalent C*-algebra (as in \cite{Sc:fine2}), we may assume that $A$ is the
closure of an increasing sequence $(A_n)$ of nuclear C*-subalgebras of $A$
satisfying the UCT and with the property that each $K_*(A_n)$ is finitely
generated. In particular the map $\Gamma_n :KK(A_n,B) \to
\mathrm{Hom}_{\Lambda}(\uk(A_n),\uk(B))$ is an isomorphism by ~\eqref{UMCT}. By
the open mapping theorem all we need to prove is that
$\mathrm{ker}(\Gamma)=\bar{0}$. The inclusion
$\mathrm{ker}(\Gamma)\supset\bar{0}$  follows from the continuity of $\Gamma$.
Conversely let $[\alpha,\gamma]\in \mathrm{ker}(\Gamma)$ with $\gamma$
absorbing. Let $\fset_n \subset A_n$ be a finite subset such that the union of
$(\fset_n)$ is dense in $A$. Since the diagram
$$ \xymatrix{
{KK(A,B)}\ar[r]\ar[d]&
{\homl}\ar[d]\\
{KK(A_n,B)}\ar[r]&{\mathrm{Hom}_{\Lambda}(\uk(A_n),\uk(B))} }$$ is commutative,
we have that $[\alpha,\gamma]=0$ when regarded as an element of $KK(A_n,B)$. By
Theorem~\ref{vanish} there is a unitary $u_n \in 1+\mathcal{K}(H_B \oplus H_B)$
such that for all $a_n \in \fset_n$
$$ \|u_n \left( \alpha(a)\oplus \gamma(a)
\right)u_n^*- \gamma(a)\oplus \gamma(a)\|<1/n.
$$
Therefore
\begin{equation*}
\lim_{n \to \infty} \|u_n \left( \alpha(a)\oplus \gamma(a)
\right)u_n^*- \gamma(a)\oplus \gamma(a)\|=0
\end{equation*}for all $a \in
A$, hence $d([\alpha,\gamma],0)=0$ and $[\alpha,\gamma] \in
\bar{0}$. \qed\end{pf}
\begin{proposition}\label{GH} Let $A$ and $B$ be separable C*-algebras and assume
that $A$ satisfies the UCT and that the group $K_*(B)$ is finitely
generated. Then for any subgroup $G$ of $KK(A,B)$ and any $\ep>0$
there is a finitely generated subgroup $H$ of $G$ which is
$\ep$-dense in $G$, i.e. for every $x \in G$ there is $y \in H$
such that $d(x,y)<\ep$.
\end{proposition}
\begin{pf}
Let $U=\{z \in KK(A,B):d(z,0)<\ep\}.$ Since the map
$$\Gamma:KK(A,B) \to \homl$$ is open, there exists an integer $m \geq 0$ and
$t_1,\dots,t_n \in \uk(A)_m$ such that $$\{\alpha\in \homl:
\alpha(t_1)=\cdots=\alpha(t_n)=0\}\subset \Gamma(U).$$ Here
$\uk(A)_m$ denotes   the subgroup of $\uk(A)$ generated by
$K_*(A;\mathbb{Z}/k)$ with $k \leq m$.
 Let
$\Gamma_n:G \to \prod_{i=1}^n \uk(B)_m$ be defined by
$\Gamma_n(x)=(\underline{x}(t_1),\dots,\underline{x}(t_n))$. Since
$K_*(B)$ is abelian and finitely generated so is $\uk(B)_m$ and its
subgroup $\Gamma_n(G)$. Therefore there is a finitely generated
subgroup $H$ of $G$ such that $\Gamma_n(G)=\Gamma_n(H)$. In
particular for any $x \in G$ there is $y \in H$ such that
$\underline{x}(t_i)=\underline{y}(t_i)$ for all $i$, $1 \leq i \leq
n$. Therefore $\underline{x}-\underline{y} \in \Gamma(U)$,  hence
$x-y \in U+\bar{0}$. We conclude that $d(x,y)=d(x-y,0)<\ep$.
\qed\end{pf} Let us recall that a C*-algebra is called nuclearly
embeddable if it has a faithful nuclear representation on a Hilbert
space. Kirchberg proved that a separable C*-algebra is nuclearly
embeddable if and only if is exact. A C*-algebra $A$ is called
residually finite dimensional (abbreviated RFD) if the finite
dimensional representations of $A$ separate the points of $A$.
 Using notation  introduced before Proposition
\ref{abs-nuc} we have:
\begin{theorem}\label{point-approx} Let $A$ be a separable unital
exact RFD C*-algebra satisfying the UCT. For any finite subset
$\fset$ of $A$ and any $\ep>0$ there are unital finite dimensional
irreducible $*$-representations $\pi_1,\dots,\pi_r$ such that for
any unital finite dimensional $*$-representation $\pi:A \to
\LLL(H_\pi)$,
 $$\pi \oplus m_1\pi_1\oplus\cdots \oplus
m_r\pi_r\SF{\ep}k_1\pi_1\oplus\cdots \oplus k_r\pi_r$$ for some
nonnegative integers $m_1,\dots m_r,k_1,\dots k_r$.\end{theorem}
\begin{pf} Let $\fdr$ denote the set of unital finite
dimensional $*$-representations of $A$. If $\pi \in\fdr$, we
denote by $[\pi]$ its class in $KK(A,\CCC)$. From the definition
of the metric $d$ we derive the following observation. Given
$\fset$ and $\ep$ as in the statement, there is $\ep_0>0$ such
that if $\pi$ and $\pi'$ are unital finite dimensional
$*$-representations of $A$ on the same space $H_\pi$ with
$d([\pi],[\pi'])<\ep_0$ then for any unitally absorbing $*$-\ho\
$\gamma:A \to\LLL(H)$ there is a unitary $u \in 1+\KKK(H_\pi\oplus
H)$ such that
$$\|\pi(a)\oplus \gamma(a)-u(\pi'(a)\oplus \gamma(a))u^*\|<\ep$$
for all $a \in \fset$. Since $A$ is separable there is a sequence
$(\pi_n)_{n=1}^\infty$ in $\fdr$ whose unitary orbit is dense in
$\fdr$ in the point-norm topology. This means that for any $\pi \in
\fdr$, any finite subset $\fset$ of $A$ and any $\ep>0$,
$\pi\SF{\ep}\pi_n$ for some $n$. Consequently it suffices to prove
the theorem only for representations $\pi$ that appear in the
sequence $(\pi_n)_{n=1}^\infty$.   We may assume that each $\pi_n$
is repeating infinitely many times. Let $G$ be the subgroup of
$KK(A,\CCC)$ generated by the set $\{[\pi_n]: n\geq 1\}$. By
Proposition~\ref{GH} there is a finitely generated subgroup $H$ of
$G$ that is $\ep_0$-dense in $G$. Therefore
    there is $r$ such that $H$
is generated by $[\pi_1],\dots,[\pi_r]$.  Fix a unitally absorbing
$*$-homomorphism $\gamma:A \to \LLL(H)$. Since $A$ is nuclearly
embeddable, by enlarging $r$, we may arrange that
\begin{equation}\label{apiax}
 \gamma \SF{\ep} \infty\cdot(\pi_1\oplus\cdots \oplus \pi_r)
\end{equation}
 by an approximation
result of \cite{Dad;apqd}; see also \cite[Prop.~6.1]{Dad:shortrfd} for a more
direct proof. Let $\pi$ be as in the statement of the theorem. We may assume
that $\pi$ appears in the sequence $(\pi_n)_{n=1}^\infty$ and therefore its
K-homology class $[\pi]$ belongs to $G$. It follows that there is $h \in H$
with $d([\pi],h)<\ep_0$. Thus there are positive integers $m_1,\dots
m_r,k_1,\dots k_r$ such that
 $$d([\pi
\oplus m_1\pi_1\oplus\cdots \oplus m_r\pi_r],[k_1\pi_1\oplus\cdots
\oplus k_r\pi_r])<\ep_0.$$ By our choice of $\ep_0$ this implies
that there is a unitary $u$ of the form $1+$\emph{compact} such
that
$$\|\pi(a)
\oplus m_1\pi_1(a)\oplus\cdots \oplus m_r\pi_r(a)\oplus
\gamma(a)-u(k_1\pi_1(a)\oplus\cdots \oplus k_r\pi_r(a)\oplus
\gamma(a))u^*\|<\ep$$ for all $a \in \fset$. Using \eqref{apiax} and
compressing by a suitable finite dimensional projection $e$ we
obtain that there exist a positive integer $N$ and a unitary $v$
close to $eue$ such that, if $M_i=m_i+N$ and $K_i=k_i+N$, then
$$ \|\pi(a) \oplus M_1\pi_1(a)\oplus\cdots \oplus
M_r\pi_r(a) -v(K_1\pi_1(a)\oplus\cdots \oplus K_r\pi_r
(a))v^*\|<3\ep$$
 This
concludes the proof. \qed\end{pf}
 If $A$ is unital, the subgroup of
$K_0(\CCC)=\mathbb{Z}$ generated by $\{[\pi(1_A)]:\pi \in \fdr\}$
is isomorphic to $d\mathbb{Z}$ for some integer $d \ge 1$. The
number $d$ is a topological invariant of $A$ and is denoted by
$d(A)$.
\begin{theorem}\label{RRFFDD} Let $A$ be a separable exact RFD C*-algebra satisfying
the UCT. Then $A$ embeds in the UHF C*-algebra of type $2^\infty$
denoted by $B$. If $A$ is unital then it embeds as a unital
C*-subalgebra in $M_{d(A)}(B)$.
\end{theorem}
\begin{pf} By adding a unit to $A$ (whether or not $A$ has
already a unit) we have $d(\tilde A)=1$. Thus it suffices to prove
only the second part of the theorem. Let $(\fset_n)_{n=1}^\infty$ be
an increasing sequence of finite subsets of $A$ whose union is dense
in $A$ and let $\ep_n=1/2^n$. By Theorem~\ref{point-approx} there
exist a sequence $(\pi_n)_{n=1}^\infty$ in $\fdr$ and integers
$0<r(1)<r(2)<\dots<r(n)<\dots$, such that if $\mathcal{R}_n\subset
\fdr$ consists of all unital representations unitarily equivalent to
representations of the form $k_1\pi_1\oplus\cdots\oplus
k_{r(n)}\pi_{r(n)}$ with $k_i>0$, then for any $\pi\in\fdr$ there
are $\alpha,\beta \in \mathcal{R}_n$ with $\pi\oplus \alpha
\SFn{\ep_n} \beta$. After changing notation if necessary, we may
assume that there is $\gamma_1\in \mathcal{R}_1$, $\gamma_1:A \to
M_{k(1)}(\CCC)$ such that $k(1)=2^m d(A)$ for some positive integer
$m$. We will construct inductively a sequence of unital $*$-\hos\
$\gamma_n:A \to M_{k(n)}(\CCC)$ with $\gamma_n \in \mathcal{R}_n$
and  such that $\|\gamma_{n+1}(a)-m(n)\gamma_n(a)\|< \ep_n$ for all
$a \in \fset_n$, where $m(n)$ is some power of $2$ and
$k(n+1)=m(n)k(n)$. Note that $\gamma_n$ will satisfy $\lim_{n\to
\infty}\|\gamma_n(a)\|=\|a\|$ for all $a \in A$ since the sequence
$(\pi_n)_{n=1}^\infty$ separates the elements of $A$. Suppose that
$\gamma_1,\dots\gamma_n$ were constructed. Pick some $\pi \in
\mathcal{R}_{n+1}$. Then $\pi\oplus\alpha\SFn{\ep_n} \beta$ for some
$\alpha,\beta \in \mathcal{R}_n$. Since $\gamma_n \in
\mathcal{R}_n$, there exists a power of $2$ denoted by $m(n)$ and
$\beta' \in \mathcal{R}_n$ such that $\beta\oplus \beta'$ is
unitarily equivalent to $m(n)\gamma_n$ hence
 $\pi\oplus\alpha\oplus\beta'\SFn{\ep_n} m(n)\gamma_n$.
 It follows that there is a finite dimensional unitary $u$ such
 that
$\|u(\pi\oplus\alpha\oplus\beta')(a)u^*-m(n)\gamma_n(a)\|<\ep_n$
for all $a \in \fset_n$.
 Setting $\gamma_{n+1}=u(\pi\oplus\alpha\oplus\beta')u^*$ we
 complete the induction process.

 Let $\iota_n:M_{k(n)}(\CCC) \hookrightarrow
 \underset{\longrightarrow}{\lim}\,\,M_{k(n)}(\CCC)\cong M_{d(A)}(B)$
 be the canonical inclusion.
 Having the sequence $\gamma_n$
available, we construct a unital embedding $\gamma:A \to
M_{d(A)}(B)$
 by defining $\gamma(a)$, $a \in \cup_{n=1}^\infty\fset_n$, to be
the limit of the Cauchy sequence $(\iota_n\gamma_n(a))$ and then
extend to $A$ by continuity. \qed\end{pf}
\begin{remark}  The AF-embeddability of
a separable nuclear RFD C*-algebra satisfying the UCT was proved in
\cite{Lin:rfd}. The approximation property given by
Theorem~\ref{point-approx} is  a stronger property than
UHF-embeddability. It is significant that it holds for exact
C*-algebras since as noted in \cite{Dad:shortrfd} the
UHF-embeddability of the cone of an exact separable RFD C*-algebra
(which satisfies the UCT by virtue of being contractible) implies
Kirchberg's fundamental characterization of exact separable
C*-algebras as subquotients of UHF algebras \cite{Kir-17}.
Subsequently Ozawa proved that AF-embeddability of separable exact
C*-algebras is a homotopy invariant \cite{Oz}.
\end{remark}
A locally compact group $G$ is called maximally almost periodic
(abbreviated MAP) if it has a separating family of finite
dimensional unitary representations. Residually finite groups are
examples of MAP groups. If $G$ is a second countable amenable
locally compact MAP group, then $C^*(G)$ is residually finite
dimensional by \cite{BLS} and satisfies the UCT by \cite{Tu:BC}.
By Theorem~\ref{RRFFDD} we have the following.

\begin{corollary}\label{groups}
The C*-algebra of a second countable amenable locally compact MAP
group $G$ is embeddable in the UHF C*-algebra of type
$2^\infty$.\end{corollary}
\begin{remark} If in addition we assume
that $G$ is discrete, then $G$ injects in the unitary group of $B$.
Note that this result is non-trivial even for the discrete
Heisenberg group $\mathbb{H}_3$, since $\mathbb{H}_3$ does not have
injective finite dimensional unitary representations. Indeed if
$\pi:\mathbb{H}_3 \to U(n)$ is an irreducible representation and
$s,t$ are generators of $\mathbb{H}_3$ such that $r=s^{-1}t^{-1}st$
generates the center of $\mathbb{H}_3$, then $\pi(r)=\lambda 1_n$,
$\lambda \in \mathbb{C},$ hence
$\lambda^n=\det(\pi(r))=\det(\pi(s^{-1}t^{-1}st))=1$.
\end{remark}

\section{From KL-equivalence to KK-equivalence}
In this section we address the question of when the Hausdorff quotient of
$KK(A,B)$ admits an algebraic description.  The following definition due to H.
Lin appears in \cite{Lin:auct}, except that the topology considered there is
the Brown-Salinas topology, which we will show to coincide with Pimsner's
topology in the next section.
 A separable C*-algebra $A$ satisfies the AUCT
if the natural map
$$\frac{KK(A,B)}{\bar{0}} \to \homl$$
is a bijection for all separable C*-algebras $B$.

 Let $KL(A,B)$
denote the quotient group $KK(A,B)/\bar{0}$. Since the Kasparov product is
continuous, it descends to an associative product $KL(A,B)\times KL(B,C)\to
KL(A,C)$. The group $KL(A,B)$ was first introduced by R{\o}rdam
\cite{Ror:cuntz} as the quotient of $KK(A,B)$ by
$\mathrm{PExt}(K_{*-1}(A),K_*(B))$ . The assumption that $A$ satisfies the UCT
was necessary in order to make $\mathrm{PExt}(K_{*-1}(A),K_*(B))$ a subgroup of
$KK(A,B)$ via the inclusion
$$\mathrm{PExt}(K_{*-1}(A),K_*(B))\hookrightarrow \mathrm{Ext}(K_{*-1}(A),K_*(B))
\hookrightarrow K(A,B).$$ In Section~\ref{UCT} we showed that if $A$ satisfies
the UCT then  $\mathrm{PExt}(K_{*-1}(A),K_*(B))$ coincides with the closure of
zero, hence the terminology is consistent.
 Two separable C*-algebras $A$ and $B$ are
\emph{KK-equivalent}, written $A\sim_{KK}B$, if there exist $\alpha
\in KK(A,B)$ and $\beta \in KK(B,A)$ such that
$$\alpha\beta=[\mathrm{id}_A]
,\quad\beta\alpha=[\mathrm{id}_B].$$
 Similarly, $A$ is  \emph{KL-equivalent} to $B,$  written
$A\sim_{KL}B$ if there exist $\alpha \in KK(A,B)$ and $\beta \in
KK(B,A)$ such that
$$\alpha\beta-[\mathrm{id}_A]\in \bar{0},\quad\beta\alpha-[\mathrm{id}_B]\in
\bar{0}.$$ Equivalently, $A\sim_{KL}B$  if and only if there exist
$\alpha \in KL(A,B)$ and $\beta \in KL(B,A)$ such that
$$\alpha\beta=[\mathrm{id}_A]
,\quad\beta\alpha=[\mathrm{id}_B].$$ by $KL(A,B)$. Note that
KL-equivalence corresponds to the notion of isomorphism in the
category with objects separable C*-algebras and morphisms from $A$
to $B$ given by  $KL(A,B)$.

 A separable  simple unital purely infinite nuclear
C*-algebra $A$ is called a Kirchberg C*-algebra
\cite[4.3.1]{Ror:ency}. One says that $A$ is in standard form if
$[1_A]=0$ in $K_0(A)$. The following result is due to H. Lin, except
that he works with the Brown-Salinas topology.
\begin{theorem}[\cite{Lin:auct}]\label{pisun/0}
Let $A$ and $B$ be unital Kirchberg C*-algebras.

(a) Let $\varphi,\psi:A \to B$ be  unital $*$-\hos. If
$[\varphi]=[\psi]$ in  $KL(A,B)$ then
$\varphi\approx_u\psi$.

(b) Assume that $A$ and $B$ are in standard form.
If $A\sim_{KL}B$
then $A$ is isomorphic to $B$.
\end{theorem}
\begin{pf} We include a new simple proof. (a) Since the constant sequence $[\psi]$
converges to $[\varphi]$ there is $y \in KK(A,C(\bn)\otimes B)$ such
that $y(n) =[\psi]$ for $n \in \mathbb{N}$ and
$y(\infty)=[\varphi]$. Since $B \cong B \otimes \mathcal{O}_\infty$
 by Kirchberg's theorem
\cite[7.2.6]{Ror:ency}, it follows by Phillips' classification theorem
\cite[Thm.~8.2.1]{Ror:ency} and by \cite[Prop.~4.1.4]{Ror:ency} that there is a
unital $*$-\ho\ $\Psi:A \to C(\bn)\otimes B$ with $y=[\Psi]$. Note that $\Psi$
is given by a family of $*$-\hos, $\Psi=(\psi_n)_{n \in \bn}$ satisfying
\begin{equation}\label{xy}\lim_{n \to
\infty}\|\psi_n(a)-\psi_\infty(a)\|=0\end{equation}
 for all $a \in A$.
Since $[\psi_n]=y(n)=[\psi]$ it follows from \cite[Thm.~8.2.1]{Ror:ency} that
$\psi_n\approx_u\psi$ for all $n \in \mathbb{N}$ and similarly
$\psi_\infty\approx_u\varphi$ since $[\psi_\infty]=y_\infty=[\varphi]$. In
combination with \eqref{xy} this gives $\varphi\approx_u\psi$. The converse
follows from Theorem~\ref{mainconv}.

(b) Let $\alpha$ and $\beta$ be as in the definition of KL-equivalence.
Applying \cite[Thm.~8.3.3]{Ror:ency} again we lift $\alpha$ and $\beta$ to
unital $*$-\hos\ $\varphi:A \to B$ and $\psi:B \to A$ such that
$[\varphi\psi]-[\mathrm{id}_B]\in \bar{0}$ and
$[\psi\varphi]-[\mathrm{id}_A]\in \bar{0}$. From part (a) we have
$\varphi\psi\approx_u\mathrm{id}_B$ and $\psi\varphi\approx_u\mathrm{id}_A$. It
follows that $A$ is isomorphic to $B$ by Elliott's intertwining argument
\cite[2.3.4]{Ror:ency}. \qed\end{pf}

\begin{corollary}\label{kl=kk}
Two separable nuclear C*-algebras are KK-equivalent if and only if they
are KL-equivalent.
\end{corollary}
\begin{pf} Any separable nuclear C*-algebra
is KK-equivalent to a unital Kirchberg algebra in standard form
\cite[Prop.~8.4.5]{Ror:ency}. We conclude the proof by applying
Theorem~\ref{pisun/0}.\qed\end{pf}

It is known that the validity of UCT for all nuclear separable C*-algebras is
equivalent to the statement that $KK(A,A)=0$ for all nuclear separable
C*-algebras $A$ with $K_*(A)=0$ (see \cite{RosSho:UCT} and
\cite[Prop.~5.3]{Ska:K-nuclear}).  The following answers an informal question
of Larry Brown and shows that if $A$ fails to satisfy the UCT then
\mbox{$KK(A,A)/\bar{0}\neq 0$}.
\begin{corollary}
Let $A$ be a separable nuclear C*-algebra. If $KK(A,A)=\bar{0}$
then $A$ satisfies the UCT and in fact $A\sim_{KK}0$.
\end{corollary}
Next we show that a nuclear separable C*-algebra satisfies the AUCT if and only
if it satisfies the UCT. This answers a question of H. Lin \cite{Lin:auct}.

\begin{theorem}\label{sol} Let $A$ be a separable nuclear C*-algebra. The following
assertions are equivalent.

 (i) $A$ satisfies the UCT.

 (ii) $A$ satisfies the AUCT.

 (iii) $A$ is KL-equivalent to a
commutative C*-algebra.

 (iv) $A$ is KK-equivalent to a
commutative C*-algebra.
\end{theorem}
\begin{pf} (i) $\Rightarrow$ (ii) follows from
Theorem~\ref{ppp}. (ii) $\Rightarrow$ (iii) Assume that $A$
satisfies the AUCT. Let $C$ be a separable commutative C*-algebra
with $K_*(C)\cong K_*(A)$. Since $C$ satisfies the UCT, there is
 $\alpha\in KK(C,A)$ such that the induced map
$\alpha_*:K_*(C)\to K_*(A)$ is a bijection. Then
$\Gamma(\alpha):\uk(C)\to\uk(A)$ is a bijection by the five lemma.
 We denote by $\dot{\alpha}$ the
image of $\alpha$ in $KL(C,A)$. For a separable C*-algebra $B$,
consider the commutative diagram
$$ \xymatrix{
{KL(A,B)}\ar[r]\ar[d]&
{\homl}\ar[d]\\
{KL(C,B)}\ar[r]&{\mathrm{Hom}_{\Lambda}(\uk(C),\uk(B))} }$$ where the vertical
maps are $x \mapsto \dot{\alpha}x$ and composition with $\Gamma(\alpha)$. The
top horizontal map is bijective by assumption and the bottom horizontal map is
bijective by Theorem~\ref{ppp}.
 Thus the map $KL(A,B)\to KL(C,B)$
 is a bijection for all separable C*-algebras $B$.
 By the usual "category theory" argument it follows
that $\dot{\alpha}$ has an inverse $\dot{\beta}\in KL(A,C)$.

 (iii) $\Rightarrow$ (iv)  follows from Corollary~\ref{kl=kk}. (iv)
$\Rightarrow$ (i) was proved in \cite{RosSho:UCT}. \qed\end{pf}

Finally let us we mention that similar methods were used to prove that if a
nuclear separable C*-algebra $A$ can be approximated by C*-subalgebras
satisfying the UCT, then  $A$ satisfies the UCT (see \cite{Dad:uct}).

\section{KK-topology versus Ext-topology}
For separable C*-algebras $A,B$, Kasparov \cite{Kas:KK} has established an
isomorphism
$$KK(A,B) \cong \mathrm{Ext}(SA,B)^{-1}.$$ These two groups come with natural
topologies, Pimsner's topology and respectively the Brown-Salinas topology. In
these section we show that Kasparov's isomorphism is a homeomorphism. The
following result and its proof is an adaptation of \cite[Thm.~3.3]{PPV1}.
\begin{proposition}\label{p-p-v} Let $A,B$ be separable C*-algebras and
let $X$ be a compact metrizable space. Then any element
$y\in\mathrm{Ext}(A,C(X)\ot B)^{-1}$ is represented by a $*$-\ho\  $\sigma:A\to
Q(C(X)\ot B\ot \KKK)$ which lifts to a completely positive contraction
$\varphi:A \to C(X)\ot M(B\ot \KKK)\subset
 M(C(X)\ot B\ot \KKK)$.
\end{proposition}
\begin{pf}
Since $y$ is an invertible extension, $y$ is represented by some  $*$-\ho\
$\tau:A\to Q(C(X)\ot B\ot \KKK)$ which lifts to a completely positive
contraction $$\phi:A \to  M(C(X)\ot B\ot \KKK)\cong \LLL(H_{C(X)\ot B})\cong
C_s(X,\LLL(H_B)).$$ By \cite[Thm.~3]{Kas:cp}, $\phi$ dilates to  a $*$-\ho\
$\rho:A\to C_s(X,\LLL(H_{B}\oplus H_{B}))$ of the form
\[\rho(a)=\begin{pmatrix} \phi(a)&\alpha (a)\\ \beta(a) & \psi(a)
\end{pmatrix},\]
such that $\alpha(a), \beta(a)\in C_s(X,\KKK(H_{ B}))$ for all $a\in A$. After
replacing $\rho$ by $\rho \oplus \gamma$ for some $(A,C(X)\ot B)$-absorbing
$*$-\ho\ $\gamma$, we may assume that $\rho$ itself is $(A,C(X)\ot
B)$-absorbing. Let $H=H'=H_{ B}$ and consider the maps $$\widetilde{\phi}:A \to
C_s(X,\LLL(H\oplus(H\oplus H')\oplus (H\oplus H')\oplus\cdots ))$$ defined by
 $$\widetilde{\phi}(a)=\phi(a)\oplus \rho(a)\oplus \rho(a)\cdots $$
 and
$$\widetilde{\rho}:A \to C_s(X,\LLL((H\oplus H')\oplus (H\oplus
H')\oplus\cdots ))$$ defined by
 $$\widetilde{\rho}(a)= \rho(a)\oplus \rho(a)\oplus\cdots .$$
Consider also the constant unitary operator $$G\in C_s(X,\LLL((H\oplus
H')\oplus (H\oplus H')\oplus\cdots,H\oplus(H\oplus H')\oplus (H\oplus
H')\oplus\cdots))$$ defined by
$$G(x)((h_1\oplus h'_1)\oplus(h_2\oplus h'_2)\oplus\cdots)=h_1\oplus(h_2\oplus h'_1)
\oplus(h_3\oplus h'_2)\oplus\cdots.$$ Let $S\in \LLL(H\oplus H\oplus\cdots)$ be
the shift operator $S(h_1\oplus h_2 \oplus \cdots)=0\oplus h_1\oplus h_2 \oplus
\cdots.$  If $U\in C_s(X,\LLL(H\oplus H', H))$ is a unitary operator, let us
define
$$\widetilde{U}\in C_s(X,\LLL((H\oplus H')\oplus (H\oplus H')\oplus\cdots, H\oplus
H\oplus\cdots))$$ by $\widetilde{U}=U\oplus U\oplus \cdots$. The following
identity was verified in the proof of \cite[Thm.~3.3]{PPV1}:
    \begin{align}\label{ppv}
    \widetilde{U}G^*\widetilde{\phi}(a)G\widetilde{U}^*&=
    \widetilde{U}\widetilde{\rho}(a)\widetilde{U}^*-
    [ U(\alpha(a)+\beta(a)U^*\oplus
    U(\alpha(a)+\beta(a))U^*\oplus\cdots]\notag\\
     &+ [U\beta(a)U^*\oplus U\beta(a)U^*
    \oplus\cdots]\circ S^*\\ &+ [U\alpha(a)U^*\oplus U\alpha(a)U^*
    \oplus\cdots]\circ S.\notag
\end{align}
Let $\gamma:A\to \LLL(H)$ be  an $(A,B)$-absorbing $*$-\ho\ and let us define
$\rho_0:A \to C_s(X,\LLL(H))$ by $\rho_0(a)(x)=\gamma(a)$ for all $a\in A$ and
$x\in X$. By Proposition~\ref{abs-nuc}, $\rho_0$ is an $(A,C(X)\ot
B)$-absorbing $*$-\ho. Since both $\rho$ and $\rho_0$ are absorbing, there is a
unitary $U\in C_s(X,\LLL(H\oplus H', H))$ such that $U\rho(a)U^*-\rho_0(a)\in
C(X,\KKK(H))$ for all $a \in A$. This shows that
$$\widetilde{U}\widetilde{\rho}(a)\widetilde{U}^*=U\rho(a)U^*\oplus
U\rho(a)U^*\oplus\cdots$$ is a norm-continuous function of $x\in X$. Since
$\alpha(a), \beta(a)\in C(X,\KKK(H)),$ and since the map $x\mapsto U(x)$ is
strictly continuous, we see that the other three terms appearing on the right
hand side of equation~\eqref{ppv} are also norm-continuous functions of $x$.
Therefore
$$\widetilde{U}G^*\widetilde{\phi}(a)G\widetilde{U}^*\in C(X,\LLL(H\oplus
H\oplus\cdots))\cong C(X)\otimes \LLL(H_B).$$ We conclude the proof by noting
   $\widetilde{U}G^*\widetilde{\phi}(\cdot)G\widetilde{U}^*$ defines the same
 element $y\in \mathrm{Ext}(A,C(X)\ot B)^{-1}$ as $\phi$.
\end{pf}
\begin{theorem}\label{A1}
Let $A$, $B$ be separable C*-algebras and let $(x_n)$ and $x_\infty$ be
elements of $\mathrm{Ext}(A,B)^{-1}$. Then $x_n\to x_\infty$ in the
Brown-Salinas topology if and only if there is $y \in \mathrm{Ext}(A,C(\bn)
\otimes B)^{-1}$ such that $y(n)=x_n$ for all $n \in \bn$.
\end{theorem}
\begin{pf}  First we prove the implication $(\Rightarrow)$. The elements of
$\mathrm{Ext}(A,B)^{-1}$ are represented by $*$-\hos\
$$\sigma:A \to Q(B \otimes \KKK)=M(B \otimes \KKK)/B \otimes \KKK$$
which admit completely positive contractive liftings $A \to M(B \otimes \KKK)$.
Such a map $\sigma$ is called liftable. Let $(\sigma_n)$, $\sigma_\infty$ be
liftable $*$-\hos\ with $x_n=[\sigma_n]$ and $x_\infty=[\sigma_\infty]$. Since
$x_n\to x_\infty$ in the Brown-Salinas topology, if $\gamma:A \to M(B \otimes
\KKK)$ is an absorbing $*$-\ho, then there is a sequence of unitaries $u_n \in
Q(B \otimes \KKK)$ liftable to unitaries in $M(B \otimes \KKK)$ such that
$$\lim_{n\to \infty}\|u_n(\sigma_n(a)\oplus\dot{\gamma}(a))u_n^*-
\sigma_\infty(a)\oplus\dot{\gamma}(a)\|=0$$ for all $a \in A$. Since
$[\sigma_n]=[u_n(\sigma_n\oplus\dot{\gamma})u_n^*]$ and
$[\sigma_\infty]=[\sigma_\infty\oplus \dot{\gamma}]$ in
$\mathrm{Ext}(A,B)^{-1}$, without loss of generality we may assume that
\begin{equation}\label{lif}\lim_{n\to \infty}\|\sigma_n(a)-
\sigma_\infty(a)\|=0\end{equation} for all $a \in A$. Define a
$*$-\ho\
$$\eta:A \to C(\bn)\otimes Q(B \otimes \KKK)\subset
 Q(C(\bn) \otimes B\otimes \KKK),$$
 by $\eta(a)(n)=\sigma_n(a)$, $n\in \bn$.
We want to show that $\eta$ is liftable.
  For $k \in \mathbb{N}$
 define
 $\eta^{(k)}:A \to C(\bn)\otimes Q(B \otimes \KKK)\subset
 Q(B \otimes C(\bn)\otimes \KKK)$ by
 $\eta^{(k)}(a)(n)=\sigma_n(a)$ if $n \leq k$ and
 $\eta^{(k)}(a)(n)=\sigma_\infty(a)$ if $n > k$.
 Note that $\eta^{(k)}$ lifts to a completely positive contraction
 $A \to C(\bn)\otimes M(B \otimes \KKK)$. Since
 $$\lim_{k\to \infty}\|\eta^{(k)}(a)-\eta(a)\|=
 \lim_{k\to \infty}\,\sup_{n>k}\|\sigma_n(a)-\sigma_\infty(a)\|=0$$
 by \eqref{lif}, it follows by a result
 of  Arveson, \cite[Thm.~6]{Arv:ext}, that $\eta$ is liftable
 and hence $y=[\eta]\in  \mathrm{Ext}(A,C(\bn) \otimes B)^{-1}$. It is clear
 that $y(n)=x_n$ for all $n \in \bn$.
 Let us prove the converse implication $(\Leftarrow)$.
By Proposition~\ref{p-p-v} every element $y \in \mathrm{Ext}(A,C(\bn) \otimes
B)^{-1}$ is represented by a $*$-homomorphism
\begin{equation}\label{gigi}
\Phi:A \to C(\bn)\otimes Q(B \otimes \KKK)\subset Q(C(\bn) \otimes B\otimes
\KKK).\end{equation} Therefore, if $(\Phi_n)_{n \in \bn}$ are the components of
$\Phi$, then
$$\lim_{n\to \infty}\|\Phi_n(a)-\Phi_\infty(a)\|=0$$
for all $a \in A$, and hence $y(n)=[\Phi_n]$ converges to
$y(\infty)=[\Phi_\infty]$ in the Brown-Salinas topology. \qed\end{pf} Let
$\beta$ be a generator of $KK^1(S\CCC,\CCC)\cong \mathbb{Z}$. The Kasparov
product
$$KK^1(S\CCC,\CCC)\otimes KK(A,B)   \to KK^1(SA,B)$$
induces a natural isomorphism
$$\chi:KK(A,B)\ni \alpha \mapsto\beta \otimes \alpha \in  KK^1(SA,B)
\cong \mathrm{Ext}(SA,B)^{-1}.$$
\begin{corollary}\label{P=BS}
Let $A$, $B$ be separable C*-algebras. The map $\chi:KK(A,B)\to
\mathrm{Ext}(SA,B)^{-1}$ is a homeomorphism, when
    $KK(A,B)$ is given
the Pimsner topology and $\mathrm{Ext}(SA,B)^{-1}$ is endowed with  the
Brown-Salinas topology.
\end{corollary}
\begin{pf} The evaluation map at $n \in \bn$ induces a commutative diagram
$$ \xymatrix{
{KK(A,C(\bn)\otimes B)}\ar[r]^{\chi}\ar[d]&
{\mathrm{Ext}(SA,C(\bn)\otimes B)^{-1}}\ar[d]\\
{KK(A,B)}\ar[r]^{\chi}&{\mathrm{Ext}(SA, B)^{-1}} }$$ Since $\chi$ is a
bijection, the result follows from Theorems~\ref{2nd} and ~\ref{p-p-v}
\qed\end{pf}

\section{Open questions}
1. Let $A$, $B$ be  separable  C*-algebras and assume that $A$ is nuclear. Is
the polish group $KK(A,B)/\bar{0}$ totally disconnected?

2.   Let $A$ be a separable nuclear C*-algebra. Fix an invariant metric for the
topology of $K^0(A)=KK(A,\CCC)$. Is it true that for any $\ep>0$ there is a
finitely generated subgroup of  $K^0(A)$ which is $\ep$-dense in $K^0(A)$?

Both questions have positive answers if one assumes that $A$
satisfies the UCT, as seen in  Section~\ref{UCT}.

\end{document}